\documentclass[12pt]{amsart}

\usepackage{scrextend} 
\usepackage{setspace}
\usepackage[top=.8in, bottom=.7in]{geometry}
\usepackage{graphicx}
\graphicspath{ {./figures/} }
\usepackage{subcaption}
\usepackage{amsmath, amsthm, amssymb, mathrsfs}
\usepackage{dsfont}
\usepackage[dvipsnames]{xcolor}
\usepackage{blindtext}
\usepackage{hyperref}
\usepackage[english]{babel}
\usepackage{latexsym}
\usepackage{enumerate}
\usepackage{mathtools}
\usepackage{appendix}
\usepackage{cleveref}  

\usepackage{todonotes}

\hypersetup{
    colorlinks=true,
    linkcolor=blue,
    citecolor=red,
    urlcolor=cyan
}

\allowdisplaybreaks  

\newtheorem{theorem}{Theorem}[section]

\title[Weak Convergence of Stochastic Linear Schr\"{o}dinger equation]{Weak Convergence of Finite Element Approximations of Stochastic Linear Schr\"{o}dinger equation driven by additive Wiener noise}

\author[M. Prasad]{Mangala Prasad}
\address{Indian Institute of Technology Kanpur, Kanpur-208016, India,}
\email{ mangalap21@iitk.ac.in 159mangala@gmail.com}

\allowdisplaybreaks
\begin{document}
	
	%
	
	
	
	\date{\today}
	
	\thanks{}
	
	\begin{abstract}
	A standard finite element method discretizes the stochastic linear Schr\"{o}dinger equation driven by additive noise in the spatial variables. The weak convergence of the resulting approximate solution is analyzed, and it is established that the weak convergence rate is twice that of the strong convergence.
	\end{abstract}

	\maketitle
	\section{Introduction and Main Result}\label{s1}
	We consider the following Stochastic linear Schr\"{o}dinger equation with additive noise,
	\begin{equation}\label{eqn1.1}
		\begin{split}
			&du+i \Delta u\,dt =dW_1 +i\, dW_2,  \quad \text{in}\quad ( 0,\infty)\times  \mathcal{O}, \\
			&u=0, \quad \text{in}\quad  (0,\infty)\times  \partial \mathcal{O},\\
			&u(0,x)=u_0(x), \quad \text{in}\quad \mathcal{O},
		\end{split}
	\end{equation}
	where $\mathcal{O} \subset{\mathbb{R}^d}$, with $d=1,2,3$ represents a polygonal domain which is bounded and convex, and  $\partial \mathcal{O}$ is the boundary of $\mathcal{O}$. $\{W_j(t)\}_{t \geq 0}$ for $j=1,2$ be  $L^2(\mathcal{O})$-valued independent Wiener processes on a filtered probability space $\big(\Omega,\mathcal{F},P,\{\mathcal{F}_t\}_{t\geq 0}\big)$ with respect to the filtration $\{\mathcal{F}_t\}_{t\geq 0}$. Furthermore, let $u_0$ be a random variable that is measurable with respect to $\mathcal{F}_0$.
    
 The strong convergence of the stochastic Schrödinger equation has been extensively studied (see \cite{finite, chuchuhong, jianbo}). However, there is comparatively little work on the weak convergence of the stochastic Schrödinger equation. To the best of our knowledge, \cite{debussche} provides a similar result, but it imposes stronger restrictions on the test function 
$\Phi$, focusing on time discretization.  In contrast, our study concerns space discretization. Similar results related to weak convergence have been studied for the heat and wave equations (see \cite{Kovacsweakh, Kovacsweakw}). The techniques used in this article have been adapted from \cite{Kovacsweakw}.
    A function $u:[0,\infty) \times \overline{\mathcal{O}} \to \mathbb{C}$ is complex-valued function. We denote 
    \[
    \begin{split}
    &u_1(t)+i\, u_2(t):=Re(u(t))+i\,Im(u(t))=u(t)\, \text{ for } t>0,\\
    &u_{0,1}+i\,u_{0,2}:=Re(u_0)+i\,Im(u_0)=u_0
     \end{split}
     \]
We write \eqref{eqn1.1} in system form as
\begin{equation*} 
    d \begin{bmatrix}
        u_1\\
        u_2
    \end{bmatrix}+\begin{bmatrix}
        0 & -\Delta\\
        \Delta & 0
    \end{bmatrix}\begin{bmatrix}
        u_1\\
        u_2
    \end{bmatrix}dt =\begin{bmatrix}
        dW_1\\
        dW_2
    \end{bmatrix} \, \text{ for } t>0, \, \begin{bmatrix}
        u_1(0)\\
        u_2(0)
    \end{bmatrix}=\begin{bmatrix}
        u_{0,1}\\
        u_{0,2}
    \end{bmatrix}.
\end{equation*}
The given system can be expressed in an abstract formulation as
\begin{equation}\label{eqn1.2}
  dX(t)=\mathbf{A}X(t)dt+dW(t)\, \text{ for }t>0,X(0)=X_0,  
\end{equation}
where \[X:=\begin{bmatrix}
  u_1 \\
  u_2  
\end{bmatrix},\, X_0:=\begin{bmatrix}
    u_{0,1}\\
    u_{0,2}
\end{bmatrix},\, \mathbf{A}=\begin{bmatrix}
    0 & -\Lambda\\
    \Lambda & 0
\end{bmatrix},\,  dW:=\begin{bmatrix}
    dW_1 \\
    dW_2
\end{bmatrix}\] and $\Lambda=-\Delta$ is the Laplacian operator.
The equation \eqref{eqn1.2} admits a unique weak  solution, which can be expressed as: 
 \begin{equation} \label{eqn1.4}
     X(t)=E(t)X_0+\int_0^tE(t-\tau)dW(\tau), \quad t \geq 0,
 \end{equation}
 where the unitary continuous group $\{E(t)\}_{\geq 0}$ is defined in Appendix \ref{apd}.
The semidiscrete approximation of \eqref{eqn1.2} involves determining $X_h(t) \in V_h \times V_h$	such that 
\begin{equation}\label{eqn1.5}
    dX_h(t)=\mathbf{A}_hX_h(t)dt +\mathbf{B}_hdW(t),\quad t>0,\, X_h(0)=X_{h,0},
\end{equation}
where 
\[\mathbf{A}_h=\begin{bmatrix}
		0 & -\Lambda _h\\
		\Lambda _h & 0
	\end{bmatrix},\, X_h=\begin{bmatrix}
		u_{h,1} \\
		u_{h,2}
	\end{bmatrix},\, X_{h,0}=\begin{bmatrix}
		u_{h,0,1}\\
		u_{h,0,2}
	\end{bmatrix}\text{ and } \mathbf{B}_h=\begin{bmatrix}
	    \mathcal{P}_h & 0\\
        0 & \mathcal{P}_h
	\end{bmatrix}.\]
  The unique mild solution to \eqref{eqn1.5} is represented as
	\begin{equation}\label{eqn1.6}
		X_h(t)=E_h(t)X_{h,0}+\int_0^t E_h(t-\tau)B_hdW(\tau),\quad t\geq 0,
	\end{equation}
 where $\{E_h(t)\}_{t \geq 0}$ is a $C_0$ semigroup on $V_h\times V_h$. The Hilbert space $H$ is defined in Appendix \ref{apd} and space $V_h$, operators $\Lambda_h\text{ and }\mathcal{P}_h$, and the semigroup $\{E_h(t)\}_{t \geq 0}$ are defined in the Appendix \ref{apdb}. Now, we present the statement of the main result of this work.
\begin{theorem}\label{thm1.1}
  Let $X$ and $X_h$ be defined as in \eqref{eqn1.4} and  \eqref{eqn1.5}, respectively. Suppose $\theta \in [0,1],$ and covariance operators $Q_i \text{ for } i=1,\,2$ satisfy 
  \[\left\|\Lambda ^{\theta /2}Q_1^{1/2}\right\|_{HS}+\left\|\Lambda ^{\theta /2}Q_2^{1/2}\right\|_{HS}<\infty.\] 
  Additionally, assume that $\Phi \in C^2_b(H,\,\mathbb{R})$. Then there is a constant $C$, depending on $\Phi,\, X_0,\,Q_i$ but not on $T,\, h$ such that for $T\geq 0$,
  \[
  \left|\mathbb{E}(\Phi(X_h(T))-\Phi(X(T)))\right|\leq Ch^{2\theta}.
  \]
\end{theorem}
Note that the assumption on covariance operators $Q_i\text{  for }i=1,2,\, \text{ given by }$ \\$\left\|\Lambda ^{\theta /2}Q_1^{1/2}\right\|_{HS}+\left\|\Lambda ^{\theta /2}Q_2^{1/2}\right\|_{HS}<\infty$, is the same as strong convergence \cite{finite}. As a result, the weak convergence rate is twice the rate of strong convergence. The proof of this theorem is in Section \ref{thmproof}. Section \ref{s2} reviews key concepts related to Hilbert-Schmidt and trace-class operators, while Subsection \ref{ss2.1} discusses the weak error in a general framework. The appendix provides definitions of the relevant spaces, finite element approximations, and the function space 
$C_b^2(H,\mathbb{R})$.
\section{Preliminaries}\label{s2}
This section provides a summary of key operator-related concepts that are essential for establishing our main result, Theorem \ref{thm1.1}.
Consider a separable Hilbert space,  $U$,  equipped with the inner product $\langle \cdot,\cdot \rangle$, and the corresponding induced norm $\|\cdot\|$ on $U$. Let $\mathcal{B}(U)$ be the set of all bounded linear operators on $U$ endowed with the standard operator norm $\|\cdot\|_{\mathcal{B}(U)}$. An operator $T_1 \in \mathcal{B}(U)$ is called Hilbert-Schmidt, if the sum 
\[
\|T_1\|_{HS}:=\sum_{j=1}^{\infty}\|T_1e_j\|^2 \, < \infty,
\]
where $\{e_j\}_{j=1}^{\infty}$ is an orthonormal basis of $U$. The Hilbert-Schmidt norm does not depend on the choice of orthonormal basis. The set of Hilbert-Schmidt operators is denoted by $\mathcal{L}_2(U)$. If $T_1\in \mathcal{L}_2(U) \text{ and } T_2 \in \mathcal{B}(U)$, then $T_1^{\star}, \, T_1T_2,\text{ and, }T_2T_1\in \mathcal{L}_2(U)$. Moreover, we have 
\[
\begin{split}
&\|T_1^{\star}\|_{HS}=\|T_1\|_{HS},\qquad \|T_1T_2\|_{HS}\leq \|T_1\|_{HS}\|T_2\|_{\mathcal{B}(U)},\\
&\|T_2T_1\|_{HS}\leq \|T_2\|_{\mathcal{B}(U)}\|T_1\|_{HS}
\end{split}
\]
Let $\mathcal{L}_1({U})$ represent the collection of trace-class operators mapping $U$ to itself. Specifically, an operator $T_1 \text{ belongs 
 to }  \mathcal{L}_1(U)$ if $T_1\in \mathcal{B}(U)$ and there exist sequences  $\{x_k\},\, \{y_k\}\subset U$ satisfying $\sum_{k=1}^{\infty}\|x_k\|\|y_k\| < \infty$ and such that 
\[
T_1y\,=\sum_{k=1}^{\infty}\langle y,y_k\rangle x_k, \quad y\in U.
\]
It is a well-established fact that $\mathcal{L}_1(U)$ forms a Banach space under the norm
\[
\|T_1\|_{Tr} \,=\inf \left\{ \sum_{k=1}^{\infty}\|x_k\|\|y_k\|:T_1y=\sum_{k=1}^{\infty}\langle y,y_k\rangle x_k\right\}.
\]
If $T_1 \in \mathcal{L}_1(U)$, then the trace of $T_1$, is given by 
\[
Tr(T_1)=\sum_{j=1}^{\infty} \langle  T_1e_j, e_j \rangle
\]
is finite, where $\{e_j\}_{j=1}^{\infty}$ is an orthonormal basis of $U$. If the orthonormal basis is changed, the trace remains unchanged. We recall well-known facts related to the trace and the trace norm which we use frequently, see \cite[Appendix C]{prato} and  \cite[Chap. 30]{Laxpd}. If  $T_1 \in \mathcal{L}_1(U) \text{ and }T_2 \in \mathcal{B}(U)$, then the both $T_1T_2 \text{ and }T_2T_1$ belong to $ \mathcal{L}_1(U)$ and 
\begin{align}
    &Tr(T_1T_2)=Tr(T_2T_1), \label{eqn2.1}\\
   & |Tr(T_1T_2)|=|Tr(T_2T_1)| \leq \|T_1\|_{Tr}\|T_2\|_{\mathcal{B}(U)}, \label{eqn2.2}\\
&\|T_1T_2\|_{Tr}\leq\|T_1\|_{Tr}\|T_2\|_{\mathcal{B}(U)}, \qquad \|T_2T_1\|_{Tr}\leq\|T_1\|_{Tr}\|T_2\|_{\mathcal{B}(U)}. \label{eqn2.3}
\end{align}
Furthermore, it $T_1 \in \mathcal{L}_1(U),$ then its adjoint $T_1^{\star} \in \mathcal{L}_1(U)$  and 
\begin{equation} \label{eqn2.4}
 Tr(T_1)=Tr(T_1^{\star}),\qquad \|T_1\|_{Tr}=\|T_1^{\star}\|_{Tr}.
\end{equation}
If both $T_1,\, T_2 \in \mathcal{L}_2(U)$, then $T_1T_2 \in \mathcal{L}_1(U)$ and 
\begin{equation}\label{eqn2.5}
    \|T_1T_2\|_{Tr}\leq \|T_1\|_{HS}\|T_2\|_{HS}.
\end{equation}
   \subsection{Weak Error Formulation} \label{ss2.1} 
  In this subsection, we present the weak error formulation within a general framework, which is used to prove Theorem \ref{thm1.1}. Let $U_1$ and $U_2$ be two real, separable Hilbert spaces. Consider  the stochastic Cauchy problem, for $t>0$
    \begin{equation}\label{eqn2.6}
      dX(t)+\mathbf{A}X(t)dt=\mathbf{B}dW(t), \text{ and}\quad X(0)=X_0.  
    \end{equation}
    Here $-\mathbf{A}$ generates a  continuous semigroup  $\{E(t)\}_{t \geq 0}$ on $U_2$, $\mathbf{B} \in \mathcal{B}(U_1,U_2)$, $\{W(t)\}_{t \geq 0}$ is a $U_1$-valued Wiener process with covariance operator $Q$ adapted to the filtration $\{\mathcal{F}_t\}_{t \geq 0}$. The initial condition $X_0$ is a  $U_2$-valued random variable that is measurable with respect to $\mathcal{F}_0$. The covariance operator $Q$ is self-adjoint and positive semidefinite. Under suitable assumptions (see \eqref{eqn2.10} below), the equation \eqref{eqn2.6} admits a unique weak solution, and it is obtained  by
    \begin{equation}\label{eqn2.7}
        X(t) \,=E(t)X_0+\int_0^tE(t-\tau)\mathbf{B}dW(\tau), \quad t\geq 0.
    \end{equation}\par
We discretize the equation \eqref{eqn2.6} and bring it into the finite-dimensional setting. To do this, consider a collection of finite-dimensional subspaces $\mathcal{S}_h$ of space $U_2$ and a collection of operators   $\mathbf{B}_h: U_1 \to \mathcal{S}_h$. This collection depends on the parameter $h$ (defined in  Appendix \ref{apdb}).  The discretization of stochastic equation \eqref{eqn2.6} is to find $X_h(t) \in \mathcal{S}_h$ such that 
    \[
    dX_h(t)+\mathbf{A}_hX_h(t)dt=\mathbf{B}_hdW(t),\quad t>0; \quad X_h(0)=X_{h,0}, 
    \]
where $-\mathbf{A}_h$ generates a continuous semigroup$\{E_h(t)\}_{t\geq 0}$ on $\mathcal{S}_h$ and  $X_{h,0}$ is a random variable that is measurable with respect to $\mathcal{F}_0$. The unique weak solution of the above Cauchy problem  is expressed as:
\begin{equation}\label{eqn2.8}
  X_h(t) \,=E_h(t)X_{h,0}+\int_{0}^tE_h(t-\tau)\mathbf{B}_hdW(\tau) , \quad t \geq 0.
\end{equation}
 If condition \eqref{eqn2.10} holds, then, according to \cite{DaPrato1992}, the equation
    \[
    \begin{split}
        &dY(t)=E(T-t)\mathbf{B}dW(t) \quad t\in (0,T],\\
        &Y(0)=E(T)X_0,
    \end{split}
    \]
  admits a unique weak solution, which is explicitly given by
    \[
    Y(t)=E(T)X_0+\int_0^t E(T-\tau)\mathbf{B}dW(\tau) \quad t\in [0,T].
    \]
   Observe that 
$
X(T)$ coincides with 
$
Y(T)$. Similar to \eqref{eqn2.8}, the discretization of $Y(t)$ is given by  
    \[
    Y_h(t)=E_h(T)X_{h,0}+\int_0^t E_h(T-\tau)\mathbf{B}_hdW(\tau) \quad t\in [0,T].
    \]
    Moreover, we see that  $X_h(T)$ (defined in \eqref{eqn2.8} ) coincides with $Y_h(T)$. Consider the following equation 
    \[
    dZ(t)=E(T-t)\mathbf{B}dW(t), \quad t\in(r,T]; \quad Z(r)=\phi.
    \]
    Here  $\phi$ is a random variable that is measurable with respect to $\mathcal{F}_{r }$. The weak solution of the above equation is represented as: 
    \begin{equation}\label{eqn2.9}
        Z(t,r,\phi)= \phi+\int_{r}^tE(T-\tau)\mathbf{B}dW(\tau),\quad t\in [r,T].
    \end{equation}
    For $\Phi \in C_b^2(U_2,\mathbb{R})$ (defined in 
 Appendix \ref{apdc}), we introduce a function $u:U_2 \times [0,T] \to \mathbb{R}$ as 
    \[
    u(x,t)=\mathbb{E}(\Phi(Z(T,t,x))),
    \]
    where $\mathbb{E}$ denotes the expectation. It follows from \eqref{eqn2.9} that the partial derivatives $u_x,u_{xx}$ of $u$  are expressed as: 
    \[
    \begin{split}
        u_x(x,t)\,=\mathbb{E}(\Phi^{\prime}(Z(T,t,x))), \quad(x,t) \in U_2\times[0,T],\\
        u_{xx}(x,t)\,=\mathbb{E}(\Phi^{\prime \prime}(Z(T,t,x))), \quad(x,t) \in U_2\times[0,T].
    \end{split}
    \]
The function $u$ satisfies (see \cite{dapratopde}) Kolmogorov's equation 
   \[
   \begin{split}
       &u_t(x,t)+\frac{1}{2}Tr(u_{xx}(x,t)E(T-t)\mathbf{B}Q\mathbf{B}^{\star}E(T-t)^{\star})=0 \quad (x,t)\in U_2 \times[0,T),\\
       &u(x,T)=\Phi(x),\quad x\in U_2.
   \end{split}
   \]
    Next, we revisit the representation formula for the weak error as given in \cite{Kovacsweakw}.
    \begin{theorem}(\cite[Theorem 3.1]{Kovacsweakw})\label{thm2.1}
        Consider the condition 
        \begin{equation}\label{eqn2.10}
            Tr\left(\int_{0}^T E(\tau)\mathbf{B}Q\mathbf{B}^{\star}E(\tau)^{\star}\right)d\tau<\infty.
        \end{equation}
        Additionally, assume that $\Phi \in C_b^2(U_2,\mathbb{R})$ (defined in Appendix\ref{apdc}), then for $T>0$, the weak error
        \[e_h(T):=\mathbb{E}\left(\Phi(X_h(T))\right)-\mathbb{E}\left(\Phi(X(T))\right)\]
         can be expressed in the following form: \begin{equation}\label{eqn2.11}
      \begin{split}
          e_h(T)=\,&\mathbb{E}(u(Y_h(0),0)-u(Y(0),0))\\
                &+\frac{1}{2}\mathbb{E}\int_0^T Tr\bigg(u_{x x}(Y_h(\tau),\tau)\\
                & \times \left[E_h(T-\tau)\mathbf{B}_h+E(T-\tau)\mathbf{B}\right]Q[E_h(T-\tau)\mathbf{B}_h-E(T-\tau)\mathbf{B}]^{\star}\bigg)d\tau\\
      \end{split}          
     \end{equation}  
     Alternatively, the weak error can also be rewritten as: 
   \begin{equation*}                 \begin{split}
          e_h(T)   =&\mathbb{E}(u(Y_h(0),0)-u(Y(0),0))\\
                &+\frac{1}{2}\mathbb{E}\int_0^T Tr\bigg(u_{x x}(Y_h(\tau),\tau)\\
                & \times [E_h(T-\tau)\mathbf{B}_h-E(T-\tau)\mathbf{B}]Q[E_h(T-\tau)\mathbf{B}_h+E(T-\tau)\mathbf{B}]^{\star}\bigg)d\tau.
             \end{split}
      \end{equation*}  
        \end{theorem}
      This result explicitly represents the weak error in the general setting in terms of the initial condition and an integral involving the second derivative of $u$, the operators $E_h(T-\tau) \text{ and } E(T-\tau)$, and covariance operator $Q$.
        \section{Proof of Main Result} \label{thmproof}
In this section, we see the proof of the Theorem \ref{thm1.1} with the help of Theorem \ref{thm2.1}.       We set $\mathcal{S}_h:=V_h \times V_h,\, \mathbf{B}_h:=\begin{bmatrix}
            \mathcal{P}_h &0\\
            0 &\mathcal{P}_h
        \end{bmatrix}$, and $X_{h,0}= (u_{h,0,1},u_{h,0,2})^T=(\mathcal{P}_hu_{0,1},\mathcal{P}_hu_{0,2})^T.$ In our case, 
        $Q$ and $\mathbf{B}$ are given by 
         \[Q=\begin{bmatrix}
            Q_1 &0\\
            0 &Q_2
        \end{bmatrix} \text{ and } \mathbf{B}=\begin{bmatrix}
            I & 0\\
            0 & I
        \end{bmatrix}\] respectively, where $Q_1\text{ and }Q_2$ are covariance operator of Wiener processes $W_1\text{ and }W_2$, respectively. Here $I$ is the identity operator on $\dot{H}^0$. The space $\dot{H}^{\gamma}$  is defined in  Appendix \ref{apd} for $ \gamma \in  \mathbb{R}$.
\begin{proof}[\bf Proof of Theorem \ref{thm1.1}]
If  there exist $\theta \in [0,1]$ such  $\|\Lambda ^{\theta /2}Q_1^{1/2}\|_{HS}+\|\Lambda ^{\theta /2}Q_2^{1/2}\|_{HS}<\infty$, then \eqref{eqn2.10} is satisfied. Let $\phi$ be a random variable that is $\mathcal{F}_t$ measurable. Then, according to \cite[Theorem 9.8]{prato}, 
    \begin{equation} \label{eqn3.1}
u(\phi,t)=\mathbb{E}\bigg(\Phi(Z(T,t,\phi))\bigg|\mathcal{F}_t\bigg),\text{ for } t\in [0,T].
    \end{equation}
From the property of conditional expectation, it follows that:
\[
\mathbb{E}\Big(u(\phi,t)\Big)=\mathbb{E}\bigg(\mathbb{E}\bigg(\Phi(Z(T,t,\phi))\bigg|\mathcal{F}_t\bigg)\bigg)=\mathbb{E}\Big(\Phi(Z(T,t,\phi))\Big),\quad \text{ for }t\in [0,T].
\]
We use an error estimate from \cite[Theorem 1.3]{finite}. For  $x \in \dot{H}^{\theta}$ we obtain 
\begin{equation}\label{eqn3.2}
\begin{split}
    \left\|\left(C_h(t)\mathcal{P}_h-C(t\right)\right)x\| \leq Ch^{\theta }\|x\|_{\theta}, \quad t\in [0,T],\quad \theta \in [0,2],\\
    \left\|\left(S_h(t)\mathcal{P}_h-S(t)\right)x\right\|\leq Ch^{\theta }\|x\|_{\theta}, \quad t\in [0,T],\quad \theta \in [0,2].
    \end{split}
\end{equation}
Let  $w=\Lambda^{\theta /2}x$. Then from the above inequality we get
\[
\begin{split}
   \|(C_h(t)\mathcal{P}_h-C(t))\Lambda^{-\theta/2 }w\| \leq Ch^{\theta }\|w\|, \quad w\in \dot{H}^{-\theta},\\
    \|(S_h(t)\mathcal{P}_h-S(t))\Lambda^{-\theta/2 }w\| \leq Ch^{\theta }\|w\|, \quad w\in \dot{H}^{-\theta}.  
\end{split}
\]
The  operators $(C_h(t)\mathcal{P}_h-C(t))\Lambda^{-\theta/2 },\, (S_h(t)\mathcal{P}_h-S(t))\Lambda^{-\theta/2 }$ are bounded on $\dot{H}^{0}$ for $\theta \geq 0$. Therefore, in terms of the operator norm, we get
\begin{equation}\label{eqn3.3}
    \begin{split}
        \|(C_h(t)\mathcal{P}_h-C(t))\Lambda^{-\theta }\|_{\mathcal{B}(\dot{H}^{0})} \leq Ch^{2\theta},\quad t\in[0,T],\quad 0\leq \theta\leq1,\\
       \|(S_h(t)\mathcal{P}_h-S(t))\Lambda^{-\theta }\|_{\mathcal{B}(\dot{H}^{0})} \leq Ch^{2\theta},\quad t\in[0,T],\quad 0\leq \theta\leq1.  
    \end{split}
\end{equation}
The function $\Phi:H \to \mathbb{R}$ is in $C^2_b(H,\mathbb{R})$. Then, by \eqref{eqn3.1} for any $y,\, \tilde{y} \in H$,
\begin{equation}\label{eqn3.4}
    (u_x(Y(t),t), y)=\mathbb{E}\Big((\Phi^{\prime}(Z(T,t,Y(t)),y)\Big|\mathcal{F}_t\Big)
\end{equation}
and 
\begin{equation}\label{eqn3.5}
 (u_{xx}(Y(t),t)y, \tilde{y})=\mathbb{E}\Big((\Phi^{\prime \prime}(Z(T,t,Y(t))y,\tilde{y})\Big|\mathcal{F}_t\Big). 
\end{equation}
For the weak error estimate, we have from \eqref{eqn3.2} and \eqref{eqn3.4} with $0\leq \theta \leq1$
\[
\begin{split}
|\mathbb{E}(u(&Y_h(0),0)-u(Y(0),0))|\\
&=\bigg|\mathbb{E}\bigg(\int_0^1\left(u_x(Y(0)+\tau(Y_h(0)-Y(0)),0),\,Y_h(0)-Y(0)\right)d\tau\bigg)\bigg|\\
&=\left|\mathbb{E}\int_0^1\mathbb{E}\left((\Phi^{\prime}(Z(T,0,Y(0)+\tau(Y_h(0)-Y(0))),\,Y_h(0)-Y(0)) \Big|\mathcal{F}_0\right)d\tau\right|\\
&\leq \sup_{x\in H}|||\Phi^{\prime}(x)|||\,\mathbb{E}|||Y_h(0)-Y(0)|||\\
&=\sup_{x\in H}|||\Phi^{\prime}(x)|||\,\mathbb{E}|||E_h(T)X_{h,0}-E(t)X_0|||\\
&=\sup_{x\in H}|||\Phi^{\prime}(x)|||\mathbb{E}\Big( \|(C_h(t)\mathcal{P}_h-C(t))u_{0,1}\|+ \|(S_h(t)\mathcal{P}_h-S(t))u_{0,2}\|\Big)\\
& \leq \sup_{x\in H}|||\Phi^{\prime}(x)|||\,Ch^{2\theta}\mathbb{E}|||X_0|||_{2\theta}.
\end{split}
\]
We simplify the integrand to bound the second term of \eqref{eqn2.11} in the error representation in Theorem \ref{thm2.1}. For simplicity, we set $s=T-\tau$ and apply \eqref{eqn2.3} to obtain
\[
\begin{split}
   &\bigg|\mathbb{E}\Big(Tr(u_{xx}(Y_h(t),t)\,(E_h(s)\mathbf{B}_h-E(s))\, Q\, [E_h(s)\,\mathbf{B}_h+E(s)]^{\star})\Big)\bigg| \\
   &=\bigg|\mathbb{E}\Big(Tr([E_h(s)\,\mathbf{B}_h+E(s)]^{\star}u_{xx}(Y_h(t),t)\,(E_h(s)\,\mathbf{B}_h-E(s))\, Q)\Big)\bigg|\\
   &\leq \left\|[E_h(s)\,\mathbf{B}_h+E(s)]^{\star}\right\|_{\mathcal{B}(H)}\, \mathbb{E}(\|u_{xx}(Y_h(t),t)\|_{\mathcal{B}(H)})\, \\
   & \times\|(E_h(s)\,\mathbf{B}_h-E(s))D^{-\theta}\|_{\mathcal{B}(H)}\, \|D^{\theta}Q\|_{Tr}
\end{split}
\]
where $D^{\theta}=\begin{bmatrix}
    \Lambda^{\theta} & 0\\
    0 & \Lambda^{\theta}
\end{bmatrix}$ and $Q=\begin{bmatrix}
    Q_1 & 0\\
    0 & Q_2
\end{bmatrix}$.
Since $E(s)$ is a unitary operator for all $s$, so 
there exists $C>0$ such that 
\[
\|[E_h(s)\mathbf{B}_h+E(s)]^{\star}\|_{\mathcal{B}(H)} \leq C,\quad \text{for all }s.
\]
For the second term, we have 
\[
\begin{split}
    \left\|u_{xx}(Y_h(t),t)\right\|_{\mathcal{B}(H)}&=\sup_{{|||y|||\leq 1}}|||u_{xx}(Y_h(t),t)y|||\\
    &=\sup_{|||y|||\leq1}|||\mathbb{E}\Big(\Phi^{\prime\prime}(Z(T,t,Y_h(t)))y\Big|\mathcal{F}_t\Big)|||\\
   & \leq \mathbb{E}\Big(\sup_{|||y|||\leq1}|||\Phi^{\prime \prime}(Z(T,t,Y_h(t)))y|||\Big|\mathcal{F}_t\Big)
\end{split}
\]
Above we have used Fatou's Lemma and \eqref{eqn3.5}. By taking the expectations of both sides, we get 
\[
\begin{split}
     \mathbb{E}\left(\|u_{xx}(Y_h(t),t)\|_{\mathcal{B}(H)}\right)&\leq \mathbb{E}\Big (\sup_{|||y|||\leq 1}|||\Phi^{\prime \prime}(Z(T,t,Y_h(t)))y|||\Big)\\
     &=\mathbb{E}\|\Phi^{\prime \prime}(Z(T,t,Y_h(t)))\|_{\mathcal{B}(H)}\\
    & \leq \sup_{y\in H}\|\Phi^{\prime \prime}(y)\|_{\mathcal{B}(H)}
\end{split}
\]
which is finite because $\Phi \in C^2_b(H,\mathbb{R})$. Using \eqref{eqn3.3}, we obtain the following for $y=(y_1,y_2)^{T} \in H$:
\[
\begin{split}
\|(E_h(s)\,\mathbf{B}_h&-E(s))D^{-\theta}\|_{\mathcal{B}(H)}\\
&=\sup_{|||y|||\leq1} |||(E_h(s)\,\mathbf{B}_h-E(s))D^{-\theta}y|||\\
&=\sup_{|||y|||\leq1}\Big(\|(C_h(s)\mathcal{P}_h-C(s))\Lambda^{-\theta}y_1 -(S_h(s)\mathcal{P}_h-S(s))\Lambda^{-\theta}y_2\|\\
&\hspace{2cm} +\|(S_h(s)\mathcal{P}_h-S(s))\Lambda^{-\theta}y_1 -(S_h(s)\mathcal{P}_h-S(s))\Lambda^{-\theta}y_2\|\Big)\\
& \leq \sup_{|||y|||\leq1}\Big(\|(C_h(s)\mathcal{P}_h-C(s))\Lambda^{-\theta}y_1\|+\|(C_h(s)\mathcal{P}_h-C(s))\Lambda^{-\theta}y_2\|\\
&\hspace{2cm} +\|(S_h(s)\mathcal{P}_h-S(s))\Lambda^{-\theta}y_1\|+\|(S_h(s)\mathcal{P}_h-S(s))\Lambda^{-\theta}y_2\|\Big)\\
& \leq 2Ch^{2\theta}\sup_{|||y|||\leq1}(\|y_1\|+\|y_2\|)\leq 2Ch^{2\theta}.
\end{split}
\]
The quantity $\|D^{\theta}Q\|_{Tr}=\|\Lambda^{\theta}Q_1\|_{Tr}+\|\Lambda^{\theta}Q_2\|_{Tr} \leq \|\Lambda^{\theta/2}Q_1^{1/2}\|_{HS}^2+\|\Lambda^{\theta/2}Q_2^{1/2}\|_{HS}^2$ is finite. For this inequality, we have used \eqref{eqn2.5}. Combining the above estimate we get
\[
\begin{split}
    \Bigg|\mathbb{E}\Bigg(\int_0^T &Tr\Big(\,u_{xx}(Y_h(\tau),\tau)\\
    & \times[E_h(T-\tau)\,\mathbf{B}_h+E(T-\tau)]Q[E_h(T-\tau)\mathbf{B}_h-E(T-\tau)]^{\star}\Big)d\tau\Bigg)\Bigg|\\
    &\leq C h^{2\theta}\sup_{y\in H}\|\Phi^{\prime \prime}\|_{\mathcal{B}(H)}(\|\Lambda^{\theta}Q_1\|_{Tr}+\|\Lambda^{\theta}Q_2\|_{Tr}).
\end{split}
\]
This concludes the proof of the Theorem \ref{thm1.1}.
\end{proof} 
\appendix
\section{Notations, Definition of Spaces and Their Norms }\label{apd}
Let $\Lambda=-\Delta$ denote the Laplace operator with domain $D(\Lambda)=H^2(\mathcal{O}) \cap H^1_0(\mathcal{O})$. Cosider  $U=L^2(\mathcal{O})$  equipped with the standard inner product $\langle \cdot,\cdot \rangle$ and norm $\|\cdot\|.$ Define   
	\begin{equation*}
		\dot{H} ^{\gamma}:=D(\Lambda ^{\gamma/2}),\hspace{1cm}\|v\|_{\gamma}:=\|\Lambda^{\gamma/2}v\|=\left(\sum_{j=1}^{\infty}\lambda^{\gamma}_j\langle v,\phi_j\rangle^2\right)^{1/2},\quad \gamma \in \mathbb{R},\quad v\in \dot{H} ^{\gamma}, 
	\end{equation*}
	where $\{(\lambda_j,\phi_j)\}_{j=1}^{\infty}$ are the eigenpairs of $\Lambda$ with orthonormal eigenvectors. Then $\dot{H} ^{\gamma} \subset \dot{H} ^{\theta} \text{ for } \gamma \geq \theta.$ It is clear that $\dot{H} ^{0}=U,\text{ } \dot{H} ^{1}={H} ^{1}_0(\mathcal{O}), \text{ } \dot{H} ^{2}= D(\Lambda)=H^2(\mathcal{O}) \cap H^1_0(\mathcal{O})$ with equivalent norms;  see \cite{Thomee}. Moreover, the dual space of $ \dot{H} ^{-\theta}$ can be identified with $( \dot{H} ^{\theta})^{\star} \text{ for } \theta >0$. We define 
	\begin{equation}
H^{\gamma}:=\dot{H}^{\gamma}\times\dot{H}^{\gamma},\hspace{1cm} |||v|||_{\gamma}^2:=\|v_1\|^2_{\gamma}+\|v_2\|^2_{\gamma},\hspace{.4cm} \gamma \in \mathbb{R}
	\end{equation}
	and set $H=H^0=\dot{H}^{0}\times\dot{H}^{0}$ with corresponding norm $|||\cdot|||=|||\cdot|||_0$ and inner product $(\cdot, \cdot)$.
    For $\gamma\in[-1,0]$, we define the operator $(\mathbf{A},D(\mathbf{A}))$ in $H^\gamma,$
	$$D(\mathbf{A})= \left\{ x\in H^\gamma : \mathbf{A}x= \begin{bmatrix}
		- \Lambda x_2\\
		\Lambda x_1
	\end{bmatrix} \in H^\gamma=\dot{H}^\gamma \times \dot{H}^\gamma
	\right\}:=H^{\gamma+2}:=\dot{H}^{\gamma+2} \times \dot{H}^{\gamma+2}$$
 \[
	\mathbf{A}:=\begin{bmatrix}
		0 & -\Lambda\\
		\Lambda & 0
	\end{bmatrix}.
	\]
	The operator $\mathbf{A}$ generates an unitary continuous group  $E(t):=e^{t\mathbf{A}} \text{ on } H^\gamma$ and it is given by
	\begin{equation}\label{eqn1.3}
		E(t)=\begin{bmatrix}
			C(t) & -S(t)\\
			S(t) & C(t)
		\end{bmatrix},\quad t\in\mathbb{R},
	\end{equation}
	where $C(t)=\cos{(t\Lambda)} \text{ and }S(t)=\sin{(t\Lambda})$ are the  cosine and sine operators. For example, using $ \{ (\lambda _j,\phi _j)\}_{j=1}^{\infty}$ the orthonormal eigenpairs of $\Lambda$, the cosine and sine operators are given as, for $ w\in \dot{H}^\gamma$ and for $t\geq 0,$
	\begin{align*}
	C(t)w=\cos{(t\Lambda)}w=\sum_{j=1}^{\infty} \cos{(t\lambda_j)}\langle w,\phi _j\rangle \phi_j,\\
	S(t)w=\sin{(t\Lambda)}w=\sum_{j=1}^{\infty} \sin{(t\lambda_j)}\langle w,\phi _j\rangle \phi_j .
	\end{align*}
    \section{Finite Element Approximations } \label{apdb}
Let $\mathcal{T}_h$ represent a regular family of triangulations of 
 the domain $\mathcal{O}$. The parameter $h$ is defined as  $$h=\max_{T\in \mathcal{T}_h}h_T,$$
 where  $h_T=\text{diam}(T)\text{ and, } T \text{ is triangle in }\mathcal{O}$. Define $V_h$ as the space of continuous, piecewise linear functions associated with $\mathcal{T}_h$ that vanish on $\partial \mathcal{O}$.  Consequently, we have $V_h \subset H^1_0(\mathcal{O})=\dot{H}^1.$ 
	
	Since $\mathcal{O}$  is assumed to be convex and polygonal, the triangulations can be precisely aligned with $\partial \mathcal{O}$. Moreover, this ensures the elliptic regularity property $\|w\|_{H^2(\mathcal{O})} \leq C\|\Lambda w\|$ for all $w \in D(\Lambda)$; see \cite{pgrisvard} for further details.
	The  orthogonal projection  $\mathcal{P}_h:\dot{H}^0 \to V_h$ is  
 presented as
	\[
	\langle \mathcal{P}_hw,\psi\rangle=\langle w, \psi \rangle \text{ for all }  \psi \in V_h.
	\]
	A discrete variant of the norm $\|\cdot\|_{\gamma}$ is defined as 
	\[
	\|w_h\|_{h,\gamma}=\|\Lambda _h^{\gamma/2}w_h\|, \hspace{.3cm} w_h\in V_h, \gamma \in \mathbb{R},
	\]
	where $\Lambda_h: V_h \to V_h$ is the discrete Laplace defined by 
	\[
	\langle \Lambda_hw_h,\psi \rangle=\langle \nabla w_h,\nabla \psi \rangle\quad \forall \psi \in V_h.
	\]
The operator $\mathbf{A}_h=\begin{bmatrix}
    0 & -\Lambda_h\\
    \Lambda_h & 0
\end{bmatrix}$ is the discrete approximation of 
 $\mathbf{A}$ and it generates a continuous    semigroup $\{E_h(t)\}_{t \geq 0}$ which is expressed as :  
 \[E_h(t)=e^{t\mathbf{A}_h}=
	\begin{bmatrix}
		C_h(t) & -S_h(t)\\
		S_h(t) & C_h(t)
	\end{bmatrix},\quad t\geq 0.
 \]
	Here, $C_h(t)=\cos{(t\Lambda _h)} \text{ and }S_h(t)=\sin{(t\Lambda _h )}$ represent the discretization of the cosine and the sine operators, respectively. 
	Similar to the infinite-dimensional setting, let  $\{(\lambda _{h,j},\phi _{h,j})\}_{j=1}^{N_h}$ be the orthonormal eigenpairs of the discrete Laplacian $\Lambda _h$, where  $N_h = dim(V_h)$. For $t\geq 0,$ we have 
	\[
    \begin{split}
	C_h(t)w_h=\cos{(t\Lambda _h)}w_h=\sum_{j=1}^{N_h}\cos{(t\lambda_{h,j})}\langle w_{h},\phi_{h,j}\rangle \phi_{h,j}, \hspace{.2cm} w_h\in V_h,\\
    S_h(t)w_h=\sin{(t\Lambda _h)}w_h=\sum_{j=1}^{N_h}\sin{(t\lambda_{h,j})}\langle w_{h},\phi_{h,j}\rangle \phi_{h,j}, \hspace{.2cm} w_h\in V_h.
    \end{split}
	\] 
\section{Differentiability }\label{apdc}

Here, we define the space 
$C^2_b(U, \mathbb{R})$, where $U$ is separable Hilbert space equipped with an inner product $\langle \cdot,\cdot\rangle$ and norm $\|\cdot\|$. A real-valued 
 function $\Phi$ belongs to  $C^2_b(U, \mathbb{R})$ if it is twice Fr\'{e}chet differentiable,  with both its first and second derivatives being continuous and bounded. Moreover, the first derivative $D\Phi(x) \text{ at }x \in U$ of $\Phi $ can be identified with an element $\Phi^{\prime}(x) \in $  such that  
\[
D\Phi(x)y=\langle\Phi^{\prime}(x),y\rangle, \text{ for all } y\in U.
\]
Similarly, the second derivative $D^2\Phi(x)$ of $\Phi$ with a linear operator $\Phi^{\prime \prime}(x) \in \mathcal{B}(U)$ such that 
\[
D^2\Phi(x)(y,z)=\langle\Phi^{\prime \prime}(x)y,z\rangle,\text{ for all } y,\, z\in U.
\]
The above results follow from the Riesz representation theorem.

Now, we define the space $ C^2(U, \mathbb{R})$. A real valued function  $\Phi$  belongs to  $ C^2(U, \mathbb{R})$ if $\Phi \in C(U, \mathbb{R}), \, \Phi^{\prime} \in C(U,U), \text{ and } \Phi^{\prime \prime } \in C(H, \mathcal{B}(U))$. We then define 
\[
C_b^2(U):=\{\Phi\in C^2(U, \mathbb{R}):\|\Phi\|_{C_b^2(U)}<\infty\}
\]
with the seminorm 
\[
\|\Phi\|_{C_b^2(U)}:=\sup_{y\in U}\|\Phi^{\prime}(y)\|+\sup_{y\in U}\|\Phi^{\prime \prime}(y)\|_{\mathcal{B}(U)}.
\]
Note that the function 
$\Phi$ is not assumed to be bounded.

	\medskip
	\noindent
	{\bf{\Large{Acknowledgement:}}} The author gratefully acknowledges the Department of Mathematics and Statistics at the Indian Institute of Technology Kanpur for providing a supportive atmosphere for research. Appreciation is also extended to the MHRD, Government of India, for the financial assistance received via the GATE fellowship.
	

\begin{thebibliography}{10}

\bibitem{finite}
Suprio Bhar, Mrinmay Biswas, and Mangala Prasad.
\newblock Finite element approximations of stochastic linear {S}chr\"odinger equation driven by additive wiener noise.
\newblock {\em arXiv preprint arXiv:2410.06006}, 2024.

\bibitem{chuchuhong}
Chuchu Chen, Jialin Hong, and Lihai Ji.
\newblock Mean-square convergence of a symplectic local discontinuous {G}alerkin method applied to stochastic linear {S}chr\"odinger equation.
\newblock {\em IMA J. Numer. Anal.}, 37(2):1041--1065, 2017.

\bibitem{jianbo}
Jianbo Cui, Jialin Hong, and Zhihui Liu.
\newblock Strong convergence rate of finite difference approximations for stochastic cubic {S}chr\"odinger equations.
\newblock {\em J. Differential Equations}, 263(7):3687--3713, 2017.

\bibitem{prato}
Giuseppe Da~Prato and Jerzy Zabczyk.
\newblock {\em Stochastic equations in infinite dimensions}, volume~44 of {\em Encyclopedia of Mathematics and its Applications}.
\newblock Cambridge University Press, Cambridge, 1992.

\bibitem{DaPrato1992}
Giuseppe Da~Prato and Jerzy Zabczyk.
\newblock {\em Stochastic equations in infinite dimensions}, volume~44 of {\em Encyclopedia of Mathematics and its Applications}.
\newblock Cambridge University Press, Cambridge, 1992.

\bibitem{dapratopde}
Giuseppe Da~Prato and Jerzy Zabczyk.
\newblock {\em Second order partial differential equations in {H}ilbert spaces}, volume 293 of {\em London Mathematical Society Lecture Note Series}.
\newblock Cambridge University Press, Cambridge, 2002.

\bibitem{debussche}
Anne de~Bouard and Arnaud Debussche.
\newblock Weak and strong order of convergence of a semidiscrete scheme for the stochastic nonlinear {S}chr\"odinger equation.
\newblock {\em Appl. Math. Optim.}, 54(3):369--399, 2006.

\bibitem{Kovacsweakh}
Matthias Geissert, Mih\'aly Kov\'acs, and Stig Larsson.
\newblock Rate of weak convergence of the finite element method for the stochastic heat equation with additive noise.
\newblock {\em BIT}, 49(2):343--356, 2009.

\bibitem{pgrisvard}
Pierre Grisvard.
\newblock {\em Elliptic problems in nonsmooth domains}, volume~69 of {\em Classics in Applied Mathematics}.
\newblock Society for Industrial and Applied Mathematics (SIAM), Philadelphia, PA, 2011.
\newblock Reprint of the 1985 original [MR0775683], With a foreword by Susanne C. Brenner.

\bibitem{Kovacsweakw}
Mih\'aly Kov\'acs, Stig Larsson, and Fredrik Lindgren.
\newblock Weak convergence of finite element approximations of linear stochastic evolution equations with additive noise.
\newblock {\em BIT}, 52(1):85--108, 2012.

\bibitem{Laxpd}
Peter~D. Lax.
\newblock {\em Functional analysis}.
\newblock Pure and Applied Mathematics (New York). Wiley-Interscience [John Wiley \& Sons], New York, 2002.

\bibitem{Thomee}
Vidar Thom\'{e}e.
\newblock {\em Galerkin finite element methods for parabolic problems}, volume~25 of {\em Springer Series in Computational Mathematics}.
\newblock Springer-Verlag, Berlin, second edition, 2006.

\end{thebibliography}

\end{document}